\documentclass[10pt]{amsart}

\usepackage{amsxtra}
\usepackage{amssymb}

\newcommand\om{\omega}

\newcommand\la{\langle}
\newcommand\ra{\rangle}

\newcommand\hh{{\mathfrak h}} 
 
\newcommand\ggo{{\mathfrak g}}

\newcommand\RR{\mathbb R}

\newcommand\ad{\operatorname{ad}}
\newcommand\Ad{\operatorname{Ad}}

\theoremstyle{plain}

\newtheorem{thm}{Theorem}[section]

\newtheorem{prop}[thm]{Proposition}

\theoremstyle{definition}
\newtheorem{defn}[thm]{Definition}

\theoremstyle{remark}

\begin{document}
\title[Invariant metrics and Hamiltonian systems]{Invariant metrics and Hamiltonian systems}

\author{Gabriela Ovando}
\address{CIEM - Facultad de Matem\'atica, Astronom\'\i a y F\'\i sica,
Universidad Nacional de C\'or\-do\-ba, Ciudad Universitaria, C\'or\-do\-ba~5000, Argentina}
\email{ovando@mate.uncor.edu}

\date{Januar 2003}

\begin{abstract} Via a non degenerate symmetric bilinear form we identify the
 coadjoint representation  with a new representation and so we induce on the 
orbits a symplectic form. By considering Hamiltonian systems on the orbits we study some 
features of them and finally find commuting functions under the 
corresponding Lie-Poisson bracket.
\end{abstract}


\maketitle

\section{Introduction}

The study of mechanical models in the last years 
had a mayor development aimed at a unified treatment of many examples known previously and at a construction of new ones \cite{STS1}. The Toda lattice is an example in which the problem is reduced to representation theory \cite{Ko1}. This was an indication that classical completely integrable systems have an intimate  relation to Lie groups \cite{R-STS} but even the simplest systems such as the open Toda lattice require an advanced technique on representation theory (\cite{Ko1}, \cite{STS}). One of the algebraic methods used in analytical Mechanics to investigate Hamiltonian systems is related to the well known Theorem of Adler-Kostant-Symes (AKS-Theorem) (see for example \cite{Ge}). By using an ad-invariant non-degenerate symmetric bilinear form it is possible to identify the coadjoint representation from the Lie group $G$ to $Aut(\ggo^{\ast})$ with the adjoint representation from $G$ to $Aut(\ggo)$ and so to induce a symplectic form on the orbits of the adjoint representation. The AKS-Theorem studies the dynamical part of the system; it describes the solutions of the Hamiltonian systems  corresponding to  ad-invariant functions. Moreover it asserts that the ad-invariant functions are in involution with respect to the corresponding Lie-Poisson bracket. This method was succesfully applied on semisimple Lie groups attached with the Killing form and the Bruhat decomposition (as  example see the generalized Toda  lattices \cite{Ko2}). 

Our aim in this work is to give a method to construct integrable systems also by using group theory. Although the results could be applied in general, it is our interest  to see them in the case of solvable Lie groups, where the AKS theorem is difficult to be applied (in fact, there is no natural ad-invariant metric in the solvable case \cite{M-R} and this fact requires itself an study).  The idea is to identify the coadjoint representation with the ``transadjoint'' representation from the Lie group $G$ into $Aut(\ggo)$ by using a non-degenerate symmetric bilinear form (not necesarly ad-invariant). In this way we induce a symplectic structure on the orbits and we study the corresponding Hamiltonian systems. As in AKS-Theorem the invariant functions play an important role, they are in involution with respect to the corresponding Lie-Poisson bracket and it is possible to describe the solution of the corresponding Hamiltonian system.  In the final section we give some applications of these ideas. 
The first example shows the convenience of having considered the analogue of AKS for arbitrary metrics on a nilpotent Lie group. Indeed the exhibited example does not admit an ad-invariant metric  but by considering another metric we apply our results and produce a Hamiltonian system on an orbit.  In the second and third examples we apply all the machinery of our results on  fixed Lie algebras. In the second example we work on a nilpotent Lie algebra and on the third one on a solvable (not nilpotent) Lie algebra. In the nilpotent case we have a Hamiltonian system corresponding to an invariant function. This system admits two invariant functions in involution (The invariants are obtained by following \cite{C-G}). Finally we find the solution of the Hamiltonian system and study the applicability of the Liouville Theorem. In the solvable case we follow a similar procedure as in the nilpotent case by considering an invariant function but in this case the corresponding Hamiltonian system is completely integrable. 

All the examples suggest a new direction of study related with the results. Although the representations are equivalent to the coadjoint one, it could be helpful to work on the Lie algebra, specially on the study of the geometric aspects of the Hamiltonian system. A  question that arises is: what kind of properties (geometrical or others) determine that the Lie group admits many functions in involution? 

\

We will assume through this work that Lie groups and  their corresponding Lie algebras are real and finite dimensional. The Lie algebra of a Lie group $G$ will be identified with the left invariant vector fields of $G$ and denoted with greek letters as usual.

\section{Preliminaries}

A Lie group $G$ acts on its Lie algebra $\ggo$ by the adjoint representation $\Ad : G \to Aut(\ggo)$ given by $Ad(g) X = dI(g)_e X$ for all $g \in G,\, X \in \ggo$ where $I(g) : G \to G$ is defined by $I(g) h = g h g^{-1}$ for all $g, \, h \in G$. 

We use de notation $\ad_X Y$ for $[X,Y]$, where $X,Y$ are elements of the Lie algebra $\ggo$.

Let $\ggo^{\ast}$ be the dual space of $\ggo$, i. e. the space of linear
real valued functionals on $\ggo$. The coadjoint representation from  $G$ into $Aut(\ggo^{\ast})$ is defined by:
$$
(\Ad^{\ast}(g) \phi) (Y) = \phi ( Ad(g^{-1}) Y) \qquad \forall g \in G, \, \, Y
\in \ggo, \,\, \phi \in \ggo^{\ast}.
$$
Under   this  action of $G$, any point $ \phi \in \ggo^{\ast}$ defines an orbit denoted by $O_{\phi} =\{ \Ad^{\ast}(g) \phi \, , \, \text{ for all } g \in G \}$. The tangent vectors to the orbit at a point $\phi  \in \ggo^{\ast}$ are the induced vector fields  defined as usual by 
$$
{X}^{\ast}_{\phi} = \frac{d}{ds}_{|_{s=0}} \Ad^{\ast}(\exp s X) \phi = - \phi \circ ad_X, \quad  \,X \in \ggo, \quad  \phi \in \ggo^{\ast},.
$$
 On each orbit $O_{\phi}$ (of dimension 2n) we define a symplectic form  $\om$ as follows:
let $X, \, Y$ be vectors in $\ggo$ and let $X^{\ast}, \, Y^{\ast}$ be 
the induced vectors on $\ggo^{\ast}$, then
\begin{equation}
 \om _{\phi} (X^{\ast}, Y^{\ast}) = \phi ( [ X , Y ] ) \qquad  X, Y \in \ggo, \quad 
\phi \in \ggo^{\ast}.
\label{ecua2}
\end{equation}

It is known that this  defines a symplectic form on each orbit of
the coadjoint representation, that is, $\omega$ is a closed
non-degenerate bilinear two-form at the orbit such that $\omega^n$ is non zero at every point of $O$. (see for example \cite{E}).

We recall that if $(M,\om)$ is a symplectic manifold and  $f : M \to \RR$ is a smooth function on $M$, then  there is exactly one vector field $X_f$ on $M$ with the property:
$$
v_p (f) = df_p (v) = \om_p (v , X_f)  \qquad \forall \, p \in M, \,\forall \, v \in T_p M
$$



The vector field $X_f$ is called the {\it Hamiltonian}  vector field 
associated to the function $f$, and $f$ is said to be the {\it Hamiltonian} (or the 
Hamiltonian function). Moreover, the {\it Hamiltonian system} for the function $f$ 
is
\begin{equation}\label{sis}
\frac{dx}{dt} = X_f (x (t))
\end{equation} 
for $x : \RR \to M$ a curve on $M$. 

\vspace{.2cm}

A function $G:M \to \RR$ is said to be a integral of motion for (\ref{sis})
 if $G$ is constant along the solutions of (\ref{sis}).

\begin{defn} A metric on a Lie group $G$ is a smooth function $\la , \ra :TG \oplus TG \to \RR$ such that for each $g \in G$, $\la \,, \ra_g : T_gG \times T_gG \to \RR$ is a non-degenerate symmetric bilinear form.
\end{defn}

We are interested in left invariant metrics, that is, those such that left multiplication by elements of the group are isometries. In this situation  to define the metric it is enough to know its definition  at the origin $T_eG$. We identify then a left invariant metric with a non-degenerate symmetric bilinear form on the Lie algebra $\ggo$.

The dual space $\ggo^{\ast}$ can be identified with $\ggo$ by using a non-degenerate symmetric bilinear form  (not necesarily ad-invariant). In fact, 
if $\langle , \rangle$ denotes such a bilinear form, to each $ \phi \in \ggo^{\ast}$
it corresponds exactly one $Z \in \ggo$ such that $\phi = \ell_Z$, where $\ell_Z(Y) = \langle Z , Y \rangle$. In this way the coadjoint representation is equivalent to a new representation that we call the ``{\it transadjoint}'' representation, which is defined by 
 \begin{equation}\label{ecu11}
 \langle \tau(g) X, Y \rangle = \langle X , Ad(g^{-1}) Y \rangle.
\end{equation} 
In particular,
$$
\langle \tau (\exp X) Y, Z \rangle = \langle Y, Ad(\exp -X) Z \rangle = \langle Y, \sum_{j=0}^{\infty} \frac{(-\ad_X)^j}{j!} Z \rangle,
$$
and so 
\begin{equation}\label{ecu22}
\tau (\exp X) Y = \sum_{j=0}^{\infty} \frac{(-\ad^t_X)^j}{j!} Z
\end{equation}
where $\ad^t_X : \ggo \to \ggo$ is the endomorphism defined by $\langle \ad_X^t Y, Z \rangle = \langle Y, \ad_X Z \rangle$.

In the special case when $G$ is  an n-step nilpotent simply connected Lie group  the exponential map $\exp : \ggo \to G$ is a diffeomorphism from the Lie algebra $\ggo$ onto $G$ and so the transadjoint representation is defined by \ref{ecu22}.
 
\vspace{.2cm}

If $U \in \ggo$, the orbit through this point is the set
$G.U =\{\tau(g) U,\,\, \forall g \in G\}$, which is diffeomorphic to the coset $G/G_U$, with  $G_U$ the  isotropy subgroup equals $G_U =\{ g \in G \,/\,\tau(g) U = U \}$.
Each $Y \in \ggo$ corresponds to a vector field $\tilde{Y}$ on the orbit, which at a point $U \in \ggo$ satisfies 
$$
\tilde{Y_U} = \frac{d}{ds}_{|_{s=0}} \tau(\exp sY) U = - \ad^{t}_Y U.
$$
The Lie algebra of $G_U$ is the set of elements $Y \in \ggo$ such that
$$ 
0 = \frac{d}{ds}_{|_{s=0}}\tau(\exp sY ) U = ad^t_Y U 
$$

\section{Hamiltonian systems}

In this section we study some Hamiltonian systems giving the principal results. To this end we begin by constructing a symplectic manifold, which will be an orbit of the transadjoint representation, denoted by $\mathcal M$.

In fact each orbit of the transadjoint representation has a natural symplectic 
structure which is induced by the canonical symplectic structure of the 
orbits of the coadjoint representation (see (\ref{ecua2})). Explicitely for $U \in \ggo$ and $Y,Z \in \ggo$ the symplectic form on $\mathcal M$ at the point $U$ is defined by
$$
\om _ U ( \tilde{Y} , \tilde{Z} ) = \la U, [ Y, Z ] \ra
$$
In fact $\omega$ is well defined. If ${\tilde{Y}}_{1} = \tilde{Y}$ then $Y_1 = Y_1 + Y - Y_1$ and $Y-Y_1$ belongs to the Lie algebra of the isotropy group at $U$, that is $\ad^t_{Y-Y_1}U =0$ and so,  $\la U, [Y, Z]\ra = \la U, [Y_1, Z]\ra$. By using a similar argument with $Z$ one proves the assertion  and it is not hard to prove that $\omega$ is symplectic and we will call it the {\it orbit symplectic structure}.

\begin{defn} Let $(M, \omega)$ be a symplectic manifold. The Poisson bracket on $C^{\infty}(M)$ related to the symplectic form $\omega$ is defined by
$$
\{f, g\} = \omega (X_f, X_g) = X_f(g) = - X_g(f)
$$
\end{defn}

\begin{defn} Smooth functions $f_1, \hdots , f_k$ on a symplectic manifold $M$  are said to be in involution or commute if $\{f_i, f_j\}=0$ for all $1 \le i,j \le k$.
\end{defn}

A subset $N$ of a manifold $M$ is said to be {\it invariant} under the flow of $X_f$ if the solution of $x'= X_f(x)$ lies in $N$ whenever $x(0) \in N$.

\begin{defn} \label{fcompin} A Hamiltonian $f$ on a symplectic manifold $M$ (of dimension 2n) is completely integrable if there exist $f_1, \hdots , f_n$ such that

(i) $\{ f, f_i\} = 0 = \{f_i, f_j\}$ for all $1 \le i,j \le n$,

(ii) $df_1, \hdots, df_n$ are linearly independent on an open subset invariant under the flow of $X_f$.
\end{defn}

In general, $m$ commuting Hamiltonians $f_1, \hdots, f_m$ on a symplectic manifold $M$ give rise to an action of $\RR^m$ on $M$. Let $(\phi_i)_t$ be the one-parameter subgroup generated by $X_{f_i}$. then
$$ (t_1, \hdots, t_m) . p = (\phi_1)_{t_1}(\phi_2)_{t_2} \hdots (\phi_m)_{t_m}$$
defines an $\RR^m$-action on $M$. Fix constants $c_1, \hdots , c_m$, since $\{f_i, f_j\}= 0$, the set $N= \{ x\in M / f_i(x) =c_i\}$ is invariant under the $\RR^m$-action. If $N$ is compact, then the $\RR^m$-action decends to a torus action on $N$. When $m = \frac{1}{2} \dim{M}$, we get the Liouville Theorem:

\begin{thm}\label{liou}({\bf Liouville}) Suppose $f$ is completely integrable on $(M^{2n}, \omega)$ and $f_1, \hdots, f_n$ are commuting Hamiltonians such that $df_1, \hdots, df_n$ are linearly independent. Assume $F = (f_1, \hdots f_n): M \to \RR^n$ is proper then $F^{-1}(c)$ is invariant under the $\RR^n$-action and decends to a torus $T^n$-action. Let $\theta_1, \hdots, \theta_n$ denote the angle coordinates on the invariant tori. Then i) $\{f_i, f_j\}=\{\theta_i, \theta_j\}= 0$ and ii) $\{f_i, \theta_j\} = c_{ij}(F)$ for some functions $c_{ij}:\RR^n \to \RR$. In particular, the flow of $X_f$ in coordinates $(f_1, \hdots , f_n,\theta_1, \hdots \theta_n)$ is linear.
\end{thm}

A Hamiltonian system corresponding to a function $H$ is said {\it completely integrable} or (Liouville integrable) if $H$ is completely integrable and there exist angle coordinates as in Theorem (\ref{liou}).

\begin{defn} If $\la , \ra$  denotes a non-degenerate symmetric bilinear form on 
$\ggo$ and $ f  :\ggo \to  \RR$ is a function on $\ggo$, we define the 
gradient of $f$ at the point $U$, $\nabla f(U)$, by
$$
\la \nabla f (U) , Y \ra = df_U (Y) \qquad \text{ for all } \,U, Y \in \ggo
$$
\end{defn}

\vspace{.2cm}

Let us consider the restriction $ H = f _{|_{\mathcal M}}$ of a function $f : \ggo \to \RR$ to an orbit $\mathcal M$. What is the Hamiltonian vector field associated to $H$? Let $U \in \mathcal{M}$, $Y \in \ggo$, then
$$
dH_U (\tilde{Y}_U) = \la \nabla f(U),\tilde{Y} \ra = \la \nabla f(U), - 
\ad^{t}_Y U \ra = \la - [Y, \nabla f(U) ], U \ra $$
and so
$$dH_U (\tilde{Y}_U) = -\omega_U(\tilde{Y}, \tilde{\nabla f}(U))
$$
On the other hand
$$
d H_U (\tilde{Y}_U) = \omega(\tilde{Y},X_H)
$$
and so we have,
\begin{equation}\label{ecu33}
X_H (U) = - \tilde{\nabla f}(U)(U) = \ad^{t}_{(\nabla f(U))} (U)
\end{equation}

where we identified $X_H$ with a  vector field in $\ggo$ also denoted by $X_H$.

\vspace{.2cm}

{\it Example} Let $\hh$ be a Lie algebra with a metric $\la \,,\, \ra$ and let $f:\hh \to \RR$ be the function
 given by $f (V) = \la Q, V \ra$ for a fixed 
point $Q \in \hh$. Then the gradient of $f$ is $\nabla f (V) = Q$. In 
fact,
$$
\la \nabla f(V) , Y \ra = df_V (Y) = \frac{d}{ds}_{|_{s=0}} f (V + s Y ) = 
\frac{d}{ds}_{|_{s=0}} \la Q, V + sY \ra = \la Q, Y \ra
$$
Consider now the orbits by the transadjoint representation with the symplectic form $\om$ above defined. By using equation (\ref{ecu33}), the Hamiltonian vector field associated to $H = f_{|_{\mathcal M}}$ 
is
$$
X_H (U) = \ad^{t}_Q (U)
$$
and the Hamiltonian system for $H$ 
$$
\left\{
\begin{array}{rcl}
x' & = & \ad^{t}_Q (x) \\
x(0) & = & P \quad \text{ for any fixed element } P \in \hh
\end{array}
\right.
$$
has the solution
$$
x(t) = \tau (\exp tQ) P.
$$
\vspace{.2cm}

In particular if the metric $\la \,, \,\ra$ is ad-invariant, then the solution is $x(t) = Ad(\exp{tQ})P$ (compare with \cite{Ge} chapter 2).

\begin{defn} A smooth function $f : \ggo \to \RR$ is said to be 
$\tau$-invariant if
$$
 f (\tau (g) U ) = f (U) \qquad \text{ for all }  g \in G, \, U \in \ggo
$$
\end{defn}


\begin{prop} \label{37} Let $f$ be a $\tau$-invariant function, then the following relations hold:

(i) $\nabla f (\tau (g) U) = Ad (g) \nabla f (U)$ for all $g \in G$, $U \in \ggo$;

(ii) $\ad^{t}_{\nabla f(U)}  U =0$ for all $U \in \ggo$.
\end{prop}
\begin{proof} (i) For the first assertion, take $g \in G$ and differentiate $f (\tau (g) U ) = f (U) $ with respect to $U$. Thus, 
$$
\frac{d}{ds}_{|_{s=0}} f(\tau (g) ( U + sV )) =  
\frac{d}{ds}_{|_{s=0}} f(U + sV)
$$
and so
$$
\begin{array}{rcl}
df_{\tau(g) U} (\tau(g) V) & = & \la \nabla f( \tau(g) U), 
\tau (g) V \ra \\
& = & \la \Ad(g^{-1}) \nabla f (\tau (g) U) , V \ra \\
\text{On the other hand, }\\
df_U (V) & = &  \la \nabla f(U), V \ra \qquad \text{ for all }  \, V 
\end{array}
$$
and so we obtain,
$$
 \Ad(g^{-1}) \nabla f(\tau (g) U) = \nabla f(U).
$$

(ii) Let $Y \in \ggo$, then by definition we have $f(\tau (\exp sY) U )= 
f(U) \, \forall\, U, Y \in \ggo$. If we differentiate at $s=0$ we obtain,
$$
\begin{array}{rcl}
0 & = & \frac{d}{ds}_{|_{s=0}} f (\tau(\exp sY) U) = df_U (\tilde{Y}_U) \\
& = & \la \nabla f(U), -\ad^{t}_Y U \ra \\
& = & \la - \ad(Y) \nabla f (U) , U \ra \\
& = & \la \ad(\nabla f (U)) Y, U \ra \\
& = & \la Y , \ad^{t}_{(\nabla f(U))} U \ra \qquad \forall \, Y \in \ggo
\end{array}
$$
and that means
$$ 
\ad ^{t} (\nabla f(U)) (U) = 0 \qquad \forall U \in \ggo.
$$
\end{proof}

\begin{prop} Let $P$ be a vector field on $M$ and let $g$ be a curve $ g : \RR \to G$ which satisfies 
$ dL_{g^{-1}} g' = P (\tau( {g(s)}^{-1}) u_0)$. 
 Then 
$$
 u(s) = \tau ( {g(s)}^{-1} ) u_0
$$
is the solution of the equation 
\begin{equation}\label{ecu44}
u' = \ad^{t}_{P(u)}(u) \qquad u(0) = u_0
\end{equation}
\end{prop}
\begin{proof} Let $y_0 \in \ggo$ and let $f$ be the function $f(s) = \la \tau(g^{-1})u_0, y_0 \ra = \la u_0, \Ad(g(s)) y_0 \ra $. Here we are 
considering the curve $y(s) = \Ad (g(s)) y_0$ which satisfies $y ' = $

\noindent$= \Ad(g) \ad_{dL _{ g^{-1}} g'} y_0 $. If we differentiate now $f(s)$, we have
$$
\begin{array}{rcl}
f'(s) & = & \la u_0 , y'(s) \ra \\
& = & \la u_0 , \Ad(g) \ad_{dL _{ g^{-1}} g'} y_0\ra \\
& = & \la \tau(g^{-1})u_0 ,  \ad_{P(u)} y \ra \\
& = & \la -\ad^{t}_{P(u)} u,  y_0 \ra \\
\end{array}
$$
and that implies (\ref{ecu44}).
\end{proof}

\begin{prop} \label{des} Suppose $\ggo$ has a non degenerate, symmetric 
bilinear form 
$\la \,, \ra $ and that $\ggo _+, \, \ggo_-$ are Lie subalgebras of $\ggo$ such that
$$
\ggo = \ggo _+ \oplus  \ggo_-
$$ 
is a direct sum of linear subspaces. Let $\ggo_+^{\perp} = \{ Y \in \ggo \, : \, 
\la X, Y \ra = 0, \, \forall X \in \ggo_+ \}$ and  $\ggo_-^{\perp} = \{ Y \in \ggo \, : \, \la
X, Y \ra = 0, \, \forall X \in \ggo_-\}$. For $X \in \ggo$, let $\ell_X \in
\ggo^{\ast}$ be defined by $\ell_X (Y) = \la X, Y \ra$. Then

i) $\ggo$ can be decomposed as a direct sum of subspaces $\ggo = \ggo_+^{\perp} \oplus
\ggo_-^{\perp}$.

ii) The maps $\ggo_+^{\perp} \to \ggo_-^{\ast}$ and $\ggo_-^{\perp} \to
\ggo_+^{\ast}$  defined by $ X \to \ell_{X_{|_{\ggo_-}}} $ and $ X \to
\ell_{X_{|_{\ggo_+}}} $ respectively  are linear isomorphisms. So the
coadjoint action of $G_-$ on $\ggo_-^{\ast}$ (respectively of $G_+$ on
$\ggo_+^{\ast}$) induces an action of $G_-$ onto $\ggo_+^{\perp}$ (resp.
 of $G_+$ onto $\ggo_-^{\perp}$), denoted by $\tilde{\tau}$.

iii) Let $\pi_{\ggo_+^{\perp}}$ be the projection of $\ggo$ onto $\ggo_+^{\perp}$
with respect to the direct sum decomposition $\ggo = \ggo_+^{\perp} + \ggo_-^{\perp}$. If $X \in \ggo_+^{\perp}$,
$ Y \in \ggo_-$ then 
$$ 
\tilde{\tau}(\exp Y) X = \pi_{\ggo_+^{\perp}} ( \tau(\exp Y) X ).
$$

iv) The infinitesimal vector field corresponding to $Y \in \ggo_-$ at the point $X \in  \ggo_+^{\perp}$ equals 
$$
\tilde{Y} (X) = \pi_{\ggo_+^{\perp}} ( -\ad^{t}_Y X ).
$$
\end{prop}
\begin{proof} i) and ii) follow from linear algebra. Let $X \in \ggo_+^{\perp}$,
$V, Y \in \ggo_-$, a direct computation gives
$$
\begin{array}{rcl}
\Ad^{\ast}(\exp Y) \ell_X (V) & = & \la X, \Ad(\exp -Y) V \ra \\
& = & \la \tau(\exp Y) X, V \ra \\
& = & \la \pi_{\ggo_+^{\perp}} (\tau(\exp Y) X), V \ra \\
& = & \ell_{\pi_{\ggo_+^{\perp}}(\tau(\exp Y) X)}(V)
\end{array}
$$
which proves iii) and iv) follows from iii).
\end{proof}

\begin{thm} \label{teo} Let $\ggo$, $\ggo_+,\, \ggo_-$ and $\la , \ra$ be as in Proposition
\ref{des}. Let $\mathcal{M} \subset \ggo_+^{\perp}$ be a $G_-$ orbit,
$\om$ the orbit symplectic structure and $f : \ggo \to \RR$ a smooth
function. Then:
\begin{enumerate}
\item The Hamiltonian vector field of $H = f_{|_{\mathcal M}}$ at $ U \in \mathcal{M}$ is 
$$
X_H (U) = - \tilde{\nabla f} (U)_- = \pi_{\ggo_+^{\perp}}( \ad^{t}_{\nabla_{-}
f(U)} U)
$$
\item If $f : \ggo \to \RR$ is $\tau$-invariant, then the Hamiltonian
equation for $H=f_{|_{\mathcal M}}$ is
$$
\frac{du}{ds} = - \ad^{t}_{\nabla f_+(u)} u = \ad^{t}_{\nabla f_-(u)}(u)
$$
\item Let $f_1, f_2 : \ggo \to \RR$ be $\tau$-invariant functions and let
$H_i = f_{i_{|_{\mathcal M}}}$,i=1,2, the restrictions respectively. Then
$$ 
\{ H_1, H_2\} = 0 
$$
\item Assume the multiplication map from $G_+ \times G_-
\to G$ defined by $(g_+ , g_-) \to g_+ g_-$ is a
diffeomorphism, and that $f : \ggo \to \RR$ is $\tau$-invariant. Then the
initial value problem
\begin{equation}
\left\{
\begin{array}{rcl}
\frac{du}{ds} & = & -\ad^{t}_{\nabla f_{+}(u)} u \\
u(0) & = & U_0
\end{array}
\right. \label{ecua1}
\end{equation}
can be solved by factorization. In fact, if we write
$$
\exp(s \nabla f(U_0)) = g_+(s) g_-(s) \in G_+ \times G_-
$$
then
$$ u(s) = \tau(\exp  s \nabla f_+(U_0)) U_0$$
 is the solution of (\ref{ecua1}).
 \end{enumerate}
\end{thm}
\begin{proof}
 
 Let $U  \in \ggo_+^{\perp}, \, Y \in \ggo_-$, then
$$
\begin{array}{rcl}
 dH_U (\tilde{Y}(U)) &  = & \la \nabla f(U), \tilde{Y} (U) 
\ra  = \la \nabla f(U) , \pi _{\ggo _+^{\perp}} (- \ad^{t}_Y (U)) \ra \\
& = & \la \nabla f_- (U), -\ad^{t}_Y (U) \ra \\
& = & \la -ad(Y) \nabla f_- (U), U \ra \\
& = & -\omega_U(\tilde{Y},\tilde{\nabla f}_-(U))
\end{array}
$$
But on the other hand
$$
d H_U (\tilde{Y}(U)) =  \omega_U(\tilde{Y}, X_H) 
$$
which implies
$$ X_H (U) = - \tilde{\nabla f}_{-}(U) (U) = \pi 
_{\ggo_+^{\perp}}(\ad^{t}_{\nabla f_-(U)} (U) ).
$$

If $f$ is $\tau$-invariant then by Proposition \ref{37} we have
$$0  = \ad^{t}_{\nabla f(U)} U = \ad^{t}_{\nabla f_+(U)} U + 
\ad^{t}_{\nabla f_-(U)} U
$$
Now
$$ X_H(U) = \pi_{\ggo_+^{\perp}}(\ad^{t}_{\nabla f_- (U)} U) = - 
\pi_{\ggo_+^{\perp}}(\ad^{t}_{\nabla f_+ (U)} U)
$$
 and we compute for $Y \in \ggo_+$
$$
\la \ad^{t}_{\nabla f_+(U)} U, Y \ra = \la U , \ad_{\nabla f_+(U)} Y \ra = 0
$$
that means that $\ad^{t}_{\nabla f_+(U)} U \in \ggo_+^{\perp}$ for all $U$, and so
$$
 X_H (U) = - \pi_{\ggo_+^{\perp}}(\ad^{t}_{\nabla f_+ (U)} U) = 
-\ad^{t}_{\nabla f_+(U)}U =  \ad^{t}_{\nabla f_- (U)} U
$$
which proves (2).

To prove (3), let $U \in \mathcal {M} \subset \ggo_+^{\perp}$. Then we have
$$
\begin{array}{rcl}
\{H_1, H_2 \}  & = & \la U , [ {\nabla f_1}_-(U) , {\nabla f_2}_-(U)] \ra 
= \la \ad^{t}_{{\nabla f_1}_-U} U , {\nabla f_2 }_-(U) \ra \\
& = & \la - \ad^{t}_{{\nabla f_1}_+(U)} U , {\nabla f_2}_- (U) \ra = \la U 
, -[{\nabla f_1}_+ (U), {\nabla f_2}_-(U)] \ra \\
& = & \la \ad^{t}_{{\nabla f_2}_-(U)}U, {\nabla f_1}_+(U) \ra \\
& = & - \la U, [ {\nabla f_2}_+ (U), {\nabla f_1}_+ (U)] \ra = 0
\end{array}
$$

It remains to prove (4). First note that $\ad^{t}_{\nabla f (U_0)}U_0 = 0 $ implies
$$
\tau(\exp s \nabla f(U_0)) U_0 \equiv  U_0
$$
Write $g(s) = \exp s \nabla f (U_0) = g_+ (s) g_-(s)$. Since
$ \tau (g_+ g_-) U_0 \equiv U_0$, then $\tau (g_+^{-1}) U_0 = \tau(g_-) 
U_0$. But
$$
dR_{g^{-1}} g' = \nabla f (U_0),
$$
where $g ' = \frac{dg}{ds}$. We compute directly to get
$$
\begin{array}{rcl}
\nabla f( \tau (g_+^{-1}) U_0 ) & = & Ad(g_+^{-1}) \nabla f(U_0) = 
dL_{g_+^{-1}} dR_{g_+} dR_{g^{-1}} g' \\ 
& = & d L_{g_+^{-1}} dR_{g_+} dR_{g_+^{-1}} d R_{g_-^{-1}}( dR_{g_-} 
g_+^{\prime} + dL_{g_+} g_-^{\prime}) \\
& = & dL_{g_+^{-1}}(g_+^{\prime}) + d R_{g_-^{-1}} (g_-^{\prime}) \in 
\ggo_+ + \ggo_-
\end{array}
$$
So $\nabla f_+(\tau ( g_+^{-1}) U_0 ) = dL_{g_+ ^{-1}}(g_+^{\prime})$. But
$$
u' = -\ad^{t}_{dL_{g_+^{-1}}(g_+^{\prime})} u = \ad^{t}_{-\nabla 
f_+(u)} u
$$
and so $u(s) = \tau(\exp s\nabla f_+(U_0))U_0$ is the solution of (\ref{ecua1}).
\end{proof}

\

As we can see the new representation $\tau$ has a direct relation with the metric. This metric is the tool that realizes the equivalence between the coadjoint representation and  the transadjoint representation. Thus, since the representations are equivalent the invariants are the same, up to equivalence, that is, they are the invariants of the coadjoint representation. However it could be helpful to work on the Lie algebra and not on the dual space. And in this context, to produce interesting Hamiltonian systems, without the assumption of having an ad-invariant metric. 
 We are interested in particular on the solvable case. To apply the theorem we need to investigate the existence of double Lie algebras, that is, a decomposition of a given Lie algebra $\ggo$ into a direct sum as vector spaces of two Lie subalgebras, and also the existence of ``enough'' functions in involution.

\section{Applications}

The topic of having ``enough'' functions in involution (that is a maximal family of commuting functions, see Definition \ref{fcompin}) is a question to solve when we try to apply  Theorem (\ref{teo}).
In the case of nilpotent Lie groups, Corwin and Greenleaf \cite{C-G} developed a method to construct the generators of the 
ring of rational $\tau$-invariant functions. We will use this method to produce the required invariant functions on each case and show the functions we obtained. (The reader could prove by simple computation that the given functions are invariant). 

The first example below shows the convenience of having considered the analogue of AKS for arbitrary metrics. Indeed the exhibited example does not admit an ad-invariant metric but by considering another metric we may apply our results and produce Hamiltonian systems on an orbit with enough functions in involution.   In the second example we work on a nilpotent Lie algebra. We are interested in a particular Hamiltonian system corresponding to an invariant function. This function is completely integrable (see Definition \ref{fcompin}). Finally we study the applicability of the Liouville theorem according to a decomposition  of the Lie algebra $\ggo$ into a vector spaces direct sum of subalgebras $\ggo_-$ and $\ggo_+$. In the third example we apply Theorem (\ref{teo}) on a four dimensional solvable Lie algebra. This Lie algebra admits an ad-invariant metric, but by taking another metric we get an equivalent Hamiltonian system for the same Hamiltonian function on a two dimensional orbit.  Furthermore the system is completely integrable (in the Liouville sense).

\

{\bf Example i)} 
If we do not require  the metric to be ad-invariant, then it is possible to construct a Hamiltonian system with sufficiently many functions in involution. This fact will be shown in this example. 

Consider the 4-step nilpotent Lie algebra $\ggo$ generated by the elements $e_1, e_2, e_3,$ $e_4, e_5$ and the Lie bracket relations given by:
$$
[e_1,e_2] = e_3 \quad [e_1, e_3] = e_4, \quad [e_5,e_1] = e_2, \quad [e_5,e_2] = e_3 \quad [e_5, e_3] = e_4, $$
and the metric such that the set $\{e_i\}$ is an orthonormal basis.

Then for  $X = \sum x_i e_i, Y = \sum y_i e_i$, the $\ad^t$-action is the following one
$$
\begin{array}{rcl} 
\ad_X^t Y & = & (x_5y_2 - x_3 y_4 - x_2y_3) e_1 + (x_1 y_3 + x_5y_3) e_2 + (x_1y_4+x_5y_4)e_3 - \\
& & - (x_1y_2+x_3y_4+x_2y_3)e_5.
\end{array}
$$
By following the method developed by Corwin and Greenleaf in \cite{C-G} we obtain three functions which generate the ring of rational $\tau$-invariant functions. In fact, at 
 a point $x = \sum x_i e_i$ their values are given by
$$
f_1(x) = x_4 \qquad \qquad  f_2(x) = x_2 x_4 - \frac 12 x_3^2 \qquad \qquad f_3(x) = x_4^2(x_1 - x_5)- x_2 x_3 x_4 + \frac13x_3^3
$$
The gradients of the functions $f_2$ and $f_3$ are respectively 
$$\nabla f_2(x) = x_4 e_2 + x_2 e_4 -x_3 e_3,$$
$$  \nabla f_3(x) = x_4^2e_1-x_3x_4 e_2 +(x^2_3 - x_2 x_4)e_3+[2(x_1-x_5)x_4-x_2x_3]e_4-x_4^2e_5 .$$

Take the following decomposition of $\ggo$ into a direct sum (as vector spaces) of subalgebras, $\ggo = \ggo_+ \oplus \ggo_-$ where
$$\ggo_+ = span\{e_1,e_2,e_3, e_4\} \qquad \qquad  \ggo_- = span\{e_5\}.$$
It is easy to see that 
$$\ggo_+^{\perp} = \ggo_- \qquad \qquad  \ggo_-^{\perp}= \ggo_+.$$
As we said $G_+$ acts on $\ggo_-^{\perp}$. The orbits at a point $y = (y_1 , y_2, y_3,y_4,0)$ have dimension zero if $y_3 =0=y_4$ and have dimension two if $y_4 \ne 0$ or $y_3 \ne 0$. So for example for the element $p= (0,0,0,1,0)$ the orbit at $p$ is the set ${\mathcal M} = \{ (u,0,v,1,0) \,/\,u,v \in \RR\}.$
 Consider now $H_3$, the restriction of $f_3$ to a symplectic manifold $\mathcal M$, which is the orbit at a point $U=\sum U^0_i e_i\in \ggo_+$ with $U_4^0 \ne 0$. For $x \in \mathcal M$ the function $H_3={f_3}_{|_{\mathcal M}}$ takes the value $H_3(x) = (x_4^0)^2x_1 -x_2x_3x_4^0 + \frac13x_3^3$, which is not constant on  the orbit. Then the Hamiltonian system for $H_3$
with initial value  $u(0)=\sum u_i^0 e_i$ is
$$
\left\{
\begin{array}{rcl}
\frac{du}{ds} &  = &  u^2_4u_2 e_1+u_4^2u_3e_2+u_4^3e_3 \\
u(0)& = & \sum u_i^0 e_i 
\end{array}
\right.
$$
This system has a solution $u(t)$ with coordinates 
$$\begin{array}{rcl}
u_1(t) & = & (u_4^0)^7 \frac{t^3}{3!} +  \frac12(u_4^0)^4u_3^0 t^2 + (u_4^0)^2 u_2^0t  + u_1^0 \\
u_2(t) & = &\frac12(u_4^0)^5{t^2} + (u_4^0)^2 u_3^0t+u_2^0\\
u_3(t) & = & (u_4^0)^3 t+u_3^0 \\
u_4(t) & = & u_4^0.
\end{array}
$$

It is not hard to prove that the set $N=\{ y \in {\mathcal M}\,:\,H_3(y) = c\}$ for a constant $c$, is not compact, hence the system is not completely integrable (see \ref{liou}).

\

{\bf Example ii)} In this example we apply the results we obtained to a 8-dimensional nilpotent Lie algebra.  This Lie algebra admits an ad-invariant metric and was considered by the author in \cite{O1} for the applicability of the AKS-Theorem. But here we are interested in the Hamiltonian system corresponding to a $\tau$-invariant function. This function is completely integrable (see Definition \ref{fcompin}). Finally we study the applicability of the Liouville theorem according to a decomposition  of the Lie algebra $\ggo$ into a vector space direct sum   of subalgebras $\ggo_-$ and $\ggo_+$.

  Consider  the simply connected Lie group $G$ corresponding to the Lie algebra $\ggo$ generated by elements $e_1, e_2, e_3, e_4, e_5, e_6, e_7, e_8$ with the Lie-bracket relations:
$$ [e_1,e_4] = - e_6 \quad [e_1,e_6]= -e_7 \quad
[e_1, e_2] = e_5 \quad [e_1, e_3] = e_2$$
$$[e_4,e_2] = -e_8 = [e_6,e_3] \quad [e_2,e_3] = e_4 \quad [e_5, e_2 ] = e_7 \quad 
[e_5, e_3 ] = -e_6$$

We have then a 5-step nilpotent Lie algebra with a metric (which is not ad-invariant) defined by:
$$
\la e_1, e_1 \ra = \la e_2, e_2 \ra = \la e_6, e_6 \ra = \la e_8, e_8 \ra =1 = -\la e_3, e_5 \ra =- \la e_4, e_7 \ra $$

The $\ad^t$-action  of $\ggo$ on $\ggo$ is given  in terms of elements  $X = \sum x_i e_i,\, Y = \sum y_i e_i,$ by:
$$
\begin{array}{rcl}
\ad^t_X( Y ) & = & (x_2y_3 - x_3y_2+x_4y_6-x_6y_4) e_1 + (x_1y_4+x_3y_8)e_6 +\\
& & + (-x_1y_3+ x_3 y_7 - x_4 y_8  - x_5 y_4)e_2 + (x_1y_6-x_2y_8)e_7 + \\
& & + (-x_2y_4-x_3y_6)e_3 + (-x_1y_2+x_2y_7+x_5y_6+x_6y_8)e_5 
\end{array}
$$


As we said, we followed the theory  developed by Corwin and Greenleaf in \cite{C-G} and obtained  for an element $X = \sum x_i e_i$ the following invariant polynomials
$$
\begin{array}{rclrcl}
P_1(X) & = &  x_4  \quad \quad P_3(X) & = &   x_4 x_7 - x_3x_8 - \frac12{x_6^2}\\

P_2(X) & = &  x_8 \quad \quad P_4(X) & = & x_1 x_8 + x_5 x_4 + x_6 x_2 + x_7 x_3
\end{array} 
$$
The gradients of $P_3$ and $P_4$ are respectively
$$
\nabla P_3(X) = -x_4 e_4 + x_8 e_5 - x_6 e_6 - x_7 e_7 - x_3 e_8,
$$
$$
\nabla P_4(X) = x_8 e_1 + x_6 e_2 -x_4 e_3 -x_3 e_4 - x_7 e_5 + x_2 e_6 - x_5 e_7 + x_1 e_8.
$$
Consider now the following subalgebras
$$ 
\ggo_+ = span\{ e_2, e_3, e_4, e_6, e_7, e_8\}, \qquad  \ggo_- = span\{e_1, e_5\}.
$$
Then the Lie algebra $\ggo$ is a vector space direct sum  of $\ggo_+$ and  $\ggo_-$,  
that is,
$$
\ggo = \ggo_+ \oplus \ggo_- = \ggo_+^{\perp} \oplus \ggo_-^{\perp}
$$
where
$$
\ggo_+^{\perp} = span\{e_1, e_3\}, \qquad \ggo_-^{\perp} = span\{e_2, e_4, e_5, e_6, e_7, e_8\}.
$$

By  Theorem  \ref{teo}, $G_+$ acts on  $g_-^{\perp}$ by the so called transadjoint representation. Moreover the $\ad^t$-action for elements $X=\sum x_i e_i \in \ggo_+$, where $x_1 = 0 = x_5$, and $ Y = \sum y_i e_i \in \ggo_-^{\perp}$, where $y_1 = 0 = y_3$, is given by:
$$
\ad^{t}(\exp X) Y  =   (x_3 y_7- x_4 y_8) e_2  + (x_2y_7  + x_6 y_8) e_5  +  x_3  y_8 e_6  - x_2 y_8 e_7.
$$
It is not difficult to see that the orbits are 4-dimensional if $y_8 \ne 0$, are 2-dimensional if $y_8 = 0$ but $y_7 \ne 0$ and are 0-dimensional if $y_8 = 0 =y_7$. Moreover two elements $ U = \sum u_i e_i$, $V = \sum v_i e_i$ are in the same orbit if and only if $u_4 = v_4$ and $u_8 = v_8$. 

Consider now the restriction of $P_3$ and $P_4$, namely $H_3$ and $H_4$ respectively, to the orbit $\mathcal M$ of the point $X^o = \sum x_i^o e_i$, where $x_1^o = 0 = x_3^o$. Then $H_3$ and $H_4$ have the following expressions at the point $X$ on the orbit
$$
H_3(X) = x_7 x_4^o - \frac{x_6^2}{2} \qquad \qquad \qquad H_4(X) = x_4^o x_5 + x_2 x_6
$$ 
which are not trivial on the orbits. 
Consider now the Hamiltonian system corresponding to $H_4$, that is, for a curve  $x : \RR \to \mathcal M$ this is
\begin{equation} 
\label{sis2}
\left\{
\begin{array}{rcl}
\ad^t_{\nabla P_{4_-}(x)} x  & = & x_7x_4 e_2 - (x_2 x_8+ x_7 x_6)e_5 + x_8 x_4 e_6 + (x_8x_6 e_7\\
x(0) & = & x^o
\end{array}
\right.
\end{equation}
For the point $x^o  \in \ggo_-^{\perp}$, $x^o = \sum x_i^oe_i$ with $x_1^o = 0 = x_3^o$ we get 
 for the system (\ref{sis2}), the curve $x(t)$ with the following coordinates
$$
\begin{array}{rcl}
x_2(t) & = & {x_4^o}^2 {x_8^o}^2\frac{t^3}{3!} + x_6^ox_8^o x_4^o \frac{t^2}{2} + x_4^ox_7^o t + x_2^o\\
x_4(t) & = & x_4^o \\
x_5(t) & = &  -[ {x_4^o}^2 {x_8^o}^3 \frac{t^4}{3!} + \frac{2}{3} x_4^o x_6^o{x_8^o}^2 t^3 + ( x_4^o x_7^ox_8^o+ \frac{1}{2}{x_6^o}^2 x_8^o)t^2 + (x_2^o {x_8^o}+ x_7^ox_6^o) t] + x_5^o  \\
x_6(t) & = &   {x_4^o}{x_8^o} t + x_6^o \\
x_7(t) & = & x_4^o {x_8^o}^2 \frac{t^2}{2}  + x_6^o x_8^o t + x_7^o \\
x_8(t)& = &  x_8^o
\end{array}
$$
It is not hard to prove that the functions $H_3$ and $H_4$ are integrals of  motion for (\ref{sis2}).
For the initial condition $x^o = \sum x_i^o e_i$ with $x_1 = 0 =x_3$ the set $\mathcal N =\{(0, x_2, 0, x_4^o, x_5, {x_6}, x_7, {x_8^o})\}$  is invariant under the flow of $H_3$ and $H_4$. Moreover $dH_3$ and $dH_4$ are linearly independent on $\mathcal N$. And so we obtain that $H_4$ is completely integrable (Definition \ref{fcompin}) (by Theorem (\ref{teo}) the functions are in involution).
However the system is not completely integrable. In fact the set
$$M_{c_1,c_2} = \{ x \in {\mathcal M} \,/\,H_3(x) = c_1,\quad H_4(x) = c_2 \}$$ 
is not compact.

\

{\bf Example iii)} In this example we have a Hamiltonian system which is completely integrable, attachet to differents metrics on the Lie algebra. This case was constructed on a solvable Lie algebra, which admits an ad-invariant metric (and was also considered on \cite{O1}. In particular the simply-connected Lie group corresponding to this Lie algebra is a model of the time space (see \cite{M-R}).

Let $\ggo$ be the solvable Lie algebra generated by elements $e_0, e_1, e_2, e_3$ with the Lie bracket relations given by:
$$
[e_3, e_1]=e_2, \quad [e_3, e_2] = -e_1, \quad [e_1, e_2] = e_0
$$
and let $\la \, , \,\ra$ be the metric that makes of the set $\{e_0, e_1, e_2, e_3\}$ an orthonormal basis. Then the $\ad^t$-action  from $\ggo$ into $\ggo$ for elements $X = \sum x_i e_i, \, Y=\sum y_ie_i$ is given by:
$$
\ad^t_XY=(x_3y_2-x_2y_0)e_1 + (x_1y_0-x_3y_1)e_2+(x_2y_1-x_1y_2)e_3$$
By simple computation one proves that the function
$$P(X) = \frac1{2}(x_1^2 + x_2^2 + 2 x_0 x_3)$$
for $X = \sum x_i e_i$ is invariant under the corresponding $\tau$-action. The gradient of $P$ is:
$$
\nabla P(X) = x_3 e_0 + x_1 e_1 + x_2 e_2 + x_0 e_3
$$
Consider now the following subalgebras
$$ 
\ggo_+ = span\{ e_0, e_1, e_2\}, \qquad  \ggo_- = span\{e_3\}.
$$
Then the Lie algebra $\ggo$ is a direct sum as vector spaces of $\ggo_+$ and  $\ggo_-$,  
that is,
$$
\ggo = \ggo_+ \oplus \ggo_- = \ggo_+^{\perp} \oplus \ggo_-^{\perp}
$$
where
$$
\ggo_+^{\perp} = span\{e_3\}, \qquad \ggo_-^{\perp} = span\{e_0, e_1, e_2\}.
$$
By  Theorem  \ref{teo}, $G_+$ acts on  $g_-^{\perp}$ by the so called transadjoint representation. Moreover the $\ad^t$-action for elements $X=\sum x_i e_i \in \ggo_+$, where $x_3 = 0 $, and $ Y = \sum y_i e_i \in \ggo_-^{\perp}$, where $y_3 = 0$, is given by:
$$
\begin{array}{rcl}
\ad^{t}_X Y  & =  &  x_1 y_0 e_2  -  x_2y_0  e_1 
\end{array}
$$
It is not difficult to see that the orbits are 2-dimensional if $y_0 \ne 0$. Consider now the restriction of $P$ to the orbit, namely for $X \in \mathcal M$, $H(X) =\frac{1}{2}(x_1^2 + x_2^2)$. Then the Hamiltonian system of $H$ for a curve $x : \RR \to \mathcal M$, $x = x_0 e_0 + x_1 e_1 + x_2 e_  2$ is
\begin{equation}\label{sis11}\left\{
\begin{array}{rcl}
\ad^t_{\nabla H_-(x)}x & = & \ad^t_{x_0e_3}(x)= x_0 x_2 e_1 - x_0 x_1 e_2\\
x(0) & = & x^o
\end{array}\right.
\end{equation}
For the point $x^o  \in \ggo_-^{\perp}$, $x^o = \sum x_i^o e_i$ with $x_3^o = 0$ we get as solution
 of the system (\ref{sis11}) the curve $x(t)$ with the following coordinates
$$
\begin{array}{rcl}
x_1(t) & = & A \sin(x_3^o t) + B \cos(x_3^o t) \\
x_2(t) & = &  A \cos(x_3^o t) - B \sin( x_3^o t)\\
x_3(t) & = &  x_3^o
\end{array}
$$
where $A$ and $B$ are constants which depend on $x^o$. 
Now the set $\mathcal N = \{X \in \mathcal M / H(X) = c\}$ is compact and non empty for $c \ge 0$.


As we said, the considered Lie algebra admits an ad-invariant metric \cite{O}, defined by:
$$
\la e_1 , e_1 \ra = \la e_2, e_2 \ra = \la e_0, e_3 \ra = 1
$$
Thus the function $P$ is the quadratic form corresponding to the ad-invariant metric and so, the gradient of $P$ is
$$ \nabla P(X) = X$$
for a point $X = \sum x_i e_i$. Consider the same decomposition of $\ggo$ as before, $\ggo = \ggo_+ + \ggo_-$ , where 
\begin{equation}
\label{deco}
\ggo _+ =  span\{ e_0, e_1, e_2\}, \qquad  \ggo_- = span\{e_3\}.
\end{equation}
Then the Lie algebra $\ggo$ is a direct sum as vector spaces of $\ggo_+$ and  $\ggo_-$,  but also of $\ggo_+^{\perp}$ and $\ggo_-^{\perp}$
where
$$
\ggo_+^{\perp} = span\{e_0\}, \qquad \ggo_-^{\perp} = span\{e_1, e_2, e_3\}.
$$
By  AKS-Theorem  $G_+$ acts on  $g_-^{\perp}$ by the so called coadjoint representation. Moreover the $\ad$-action for elements $X=\sum x_i e_i \in \ggo_+$, where $x_3 = 0 $, and $ Y = \sum y_i e_i \in \ggo_-^{\perp}$, where $y_0 = 0$, is given by:
$$
\begin{array}{rcl}
\ad^{\ast}_X Y  & =  &  -x_1 y_3 e_2  +  x_2y_3  e_1 
\end{array}
$$
It is not difficult to see that the orbits are 2-dimensional if $y_3 \ne 0$. Consider now the restriction of $P$ to the orbit, namely for $X \in \mathcal M$, $H(X) =\frac{1}{2}(x_1^2 + x_2^2)$. Then the Hamiltonian system of $H$ for a curve $x : \RR \to \mathcal M$, $x = x_1 e_1 + x_2 e_2 + x_3 e_  3$ is
\begin{equation}\label{sis22}\left\{
\begin{array}{rcl}
\ad^{\ast}_{\nabla H_-(x)}x & = & \ad^{\ast}_{x_3e_3}(x)= -x_3 x_2 e_1 + x_3 x_1 e_2\\
x(0) & = & x^o
\end{array}\right.
\end{equation}
For the point $x^o  \in \ggo_-^{\perp}$, $x^o = \sum x_i^o$ with $x_0^o = 0$ we get 
 for the system (\ref{sis22}), the curve $x(t)$ with the following coordinates
$$
\begin{array}{rcl}
x_1(t) & = &  A \cos(x_3^o t) - B \sin( x_3^o t)\\
x_2(t) & = &  A \sin(x_3^o t) + B \cos(x_3^o t)\\
x_3(t) & = &  x_3^o
\end{array}
$$
where $A$ and $B$ are constants which depend on $x^o$. As before, the set $\mathcal N = \{X \in \mathcal M / H(X) = c\}$ is compact and non empty for $c \ge 0$.

\vspace{.2cm}

The system (\ref{sis22}) has the following matricial realization: take $M$ 
and $L$ be the following matrices
$$
M = \left(
\begin{matrix}
0 & -x_3 & 0\\
x_3 & 0 & 0 \\
0 & 0 & 0 
\end{matrix} 
\right) \qquad \quad
L = \left(
\begin{matrix}
0 & -x_3 & x_1\\
x_3 & 0 & x_2 \\
0 & 0 & 0 
\end{matrix} 
\right)
$$
Then the system can be written as $L^{\prime} = [M, L]= ML - LM$.

\vspace{.2cm}

Consider the same Lie algebra but with the metric defined by
$$
\la e_1 , e_1 \ra = \la e_2, e_2 \ra = - \la e_0, e_3 \ra = 1
$$
which is not ad-invariant and take the same decomposition of $\ggo$ into a direct sum as vector spaces $\ggo = \ggo_+ + \ggo_-$, where, $\ggo _+$ and $\ggo_-$ are as above (\ref{deco}). By considering the corresponding $\tau$-action of $G_+$ into $\ggo_-^{\perp}$ one constructs also a Hamiltonian system for the function $H$ which is the restriction of the function $P$ to the orbit. The system one obtains is equivalent to (\ref{sis11}). 

Consider the same Lie algebra $\ggo$ attached with the metric defined by
$$
- \la e_1 , e_1 \ra = \la e_2, e_2 \ra = \la e_0, e_3 \ra = 1
$$
which is not ad-invariant and take the same decomposition of $\ggo$ into a direct sum as vector spaces $\ggo = \ggo_+ + \ggo_-$, where, $\ggo _+$ and $\ggo_-$ are as above (\ref{deco}). By taking the corresponding $\tau$-action of $G_+$ into $\ggo_-^{\perp}$ one construct also a Hamiltonian system for the function  $H$ which is the restriction of the function $P$ to the orbit. The system one obtains is equivalent to (\ref{sis22}). 

\

{\it Remarks}

i) The reader could recognize that in the third example the systems are expressions of the basic equation of the theory of oscillations (see \cite{A}). 

ii) The second case in the Example iii) can be generalized to higher dimensions. In fact, let $M$ an $L$ be the following matrices on $M(2n+1, \RR)$:
$$
M = \left(
\begin{matrix}
0 & -x_3 & 0& 0 &&& &0\\
x_3 & 0 & 0 & 0 &&&&0 \\
0 & 0 & 0 & -x_3 &&&&0\\
0& 0 & x_3 & 0&&&&0\\
& & & & \ddots & & &\vdots\\
& & & & & 0 & -x_3 &0\\
& & & & 0& x_3 & 0 &0\\
0& 0& \hdots & & & && 0
\end{matrix} 
\right)
$$

$$
L = \left(
\begin{matrix}
0 & -x_3 & 0& 0 & & & &x_1\\
x_3 & 0 & 0 & 0  & & & & y_1\\
0 & 0 & 0 & -x_3 & & & & x_2\\
0& 0 & x_3 & 0 & & & & y_2\\
& & & & \ddots & &  & \vdots\\
& & & & & 0 & -x_3 & x_n\\
& & & & & x_3 & 0 & y_n\\
0& 0 & \hdots & & & & & 0
\end{matrix} 
\right)
$$
Then the equation $L^{\prime} = [M, L]= ML - LM$ give for the Hamiltonian function $H=\frac12\sum (x_i^2 + y_i^2)$ the following Hamiltonian system:

\begin{equation}\label{sis23}
\left\{
\begin{array}{rcl}
x_i'& = & -x_3 y_i\\
y_i'& = & x_3x_i\\
x_3'& =& 0 \\
x(0) & = & x^o
\end{array}\right.
\end{equation}

The system can be also realized in terms of the AKS-Theorem, by taking $\ggo$ the solvable Lie algebra generated by $X_0, X_i, Y_i, X_{n+1},\, i=1, \hdots n$, with the Lie bracket relations:
$$
[X_{n+1}, X_i]=Y_i, \quad [X_{n+1}, Y_i] = -X_i, \quad [X_i, Y_i] = X_0 \qquad \text{ for all }i
$$
and the ad-invariant metric defined for the point $X =x_0X_0 + \sum (x_iX_I+ y_i Y_i) + x_{n+1}X_{n+1}$ by the quadratic
$$q(X)= \sum (x_i^2+y_i^2)+ 2x_0x_{n+1}.$$

By taking the subalgebras 
$$\ggo_-=span\{X_0, X_i, Y_i\}_{i=1, \hdots, n} \quad \qquad \ggo_-=span\{X_{n+1}\}$$

and following a similar procedure as used in all previous examples, we get the phase space on $\ggo_-^{\perp}$. In this way we obtain a Hamiltonian system for the restriction of the function $\frac12 q(X)$ to the orbit, which is the function $H$ considered above.

This example is also known as the Hamiltonian flow on the sphere. Furthermore every quadratic, positive definite Hamiltonian function is symplectically equivalent to ones of the form $H = \sum_{j=1}^n a_j(x_j^2+y_j^2)$, which is a generalization of the considered in our example (see \cite{Mc-S}).

 \

{\bf Acknowledgement} The author thanks D. Alekseevski for helpful suggestions on the first version of this work and thanks I. Dotti Miatello for her valuable comments.

\end{document}